\documentclass{article}

\usepackage{amsthm}
\usepackage{amssymb}
\usepackage{amsmath}
\usepackage{cite}

\newtheorem{theorem}{Theorem}[section]
\newtheorem{corollary}[theorem]{Corollary}

\newtheorem{lemma}[theorem]{Lemma}
\theoremstyle{remark}
\newtheorem{remark}{Remark}[section]

\theoremstyle{definition}

\theoremstyle{definition}
\newtheorem{example}{Example}[section]

\begin{document}

\title{Diamond-$\alpha$ Jensen's Inequality
on Time Scales\thanks{This is a preprint of an article whose final and definitive form will appear in the \emph{Journal of Inequalities and Applications}, URL: {\tt http://www.hindawi.com/journals/jia/}. Accepted 07/April/2008.}}

\author{Moulay Rchid Sidi Ammi\\
\texttt{sidiammi@ua.pt}
\and Rui A. C. Ferreira\\
\texttt{ruiacferreira@yahoo.com}
\and Delfim F. M. Torres\\
\texttt{delfim@ua.pt}}

\date{Department of Mathematics\\
University of Aveiro\\
3810-193 Aveiro, Portugal}

\maketitle

%%%%%%%%%%%%%%%%%%%%%

\begin{abstract}
The theory and applications of dynamic derivatives on time scales
has recently received considerable attention. The primary purpose
of this paper is to give basic properties of diamond-$\alpha$
derivatives which are a linear combination of delta and nabla
dynamic derivatives on time scales. We prove a generalized version
of Jensen's inequality on time scales via the diamond-$\alpha$
integral and present some corollaries, including H\"{o}lder's and
Minkowski's diamond-$\alpha$ integral inequalities.
\end{abstract}

%%%%%%%%%%%%%%%%%%%%%%%%%%%%%%%%%%%%%%%%%%%%%%%%%%%%%%%%%%%%%%%%%%%%

\smallskip

\textbf{Mathematics Subject Classification 2000:} 26D15, 39A13.

\smallskip

%%%%%%%%%%%%%%%%%%%%%

\smallskip

\textbf{Keywords:} time scales, Jensen's inequality, convex
functions, dynamic derivatives.

\medskip

%%%%%%%%%%%%%%%%%%%%%%

\section{Introduction}

Jensen's inequality is of great interest in the theory of
differential and difference equations, as well as other areas of
mathematics. The original Jensen's inequality can be stated as
follows:

\begin{theorem}\emph{\cite{mit}}
\label{thm1} If $g \in C([a, b], (c, d))$ and $f \in C((c, d),
\mathbb{R} )$ is convex, then
$$
f \left( \frac{\int_{a}^{b}g(s) ds}{b-a}\right ) \leq
\frac{\int_{a}^{b}f(g(s)) ds}{b-a}.
$$
\end{theorem}

Jensen's inequality on time scales via $\Delta$-integral has
been recently obtained by Agarwal, Bohner and Peterson.

\begin{theorem}\emph{\cite{abp}}
\label{thm2} If $ g \in C_{rd}([a, b], (c, d))$ and $f \in C((c,
d), \mathbb{R}) $ is convex, then
$$
f \left( \frac{\int_{a}^{b}g(s) \Delta s}{b-a}\right ) \leq
\frac{\int_{a}^{b}f(g(s)) \Delta s}{b-a}.
$$
\end{theorem}
Under similar hypotheses, we may replace the $\Delta$-integral by
the $\nabla$-integral and get a completely analogous result
\cite{rev1:r}. The aim of this paper is to extend Jensen's
inequality to an arbitrary time scale via the diamond-$\alpha$
integral \cite{sfhd}.

There has been recent developments of the theory and applications
of dynamic derivatives on time scales. From the theoretical point
of view, the study provide a unification and extension of
traditional differential and difference equations. Moreover, it is
a crucial tool in many computational and numerical applications.
Based on the well-known $\Delta$ (delta) and $\nabla$ (nabla)
dynamic derivatives, a combined dynamic derivative, so called
$\diamondsuit_{\alpha}$ (diamond-$\alpha$) dynamic derivative, was
introduced as a linear combination of $\Delta$ and $\nabla$
dynamic derivatives on time scales \cite{sfhd}. The
diamond-$\alpha$ derivative reduces to the $\Delta$ derivative for
$\alpha =1$ and to the $\nabla$ derivative for $\alpha =0$. On the
other hand, it represents a ``weighted dynamic derivative'' on any
uniformly discrete time scale when $\alpha =\frac{1}{2}$. We refer
the reader to \cite{Stef,sQ,sfhd} for an account of the calculus
associated with the diamond-$\alpha$ dynamic derivatives.

The paper is organized as follows. In Section~\ref{sec:pre} we
briefly give the basic definitions and theorems of time scales as
introduced in Hilger's thesis \cite{h1} (see also \cite{h2,h3}).
In Section~\ref{sec:mr} we present our main results, which are
generalizations of Jensen's inequality on time scales. Some
examples and applications are given in Section~\ref{sec:app}.

%%%%%%%%%%%%%%%%%%%%%%

\section{Preliminaries}
\label{sec:pre}

A time scale $\mathbb{T}$ is an arbitrary nonempty closed subset
of real numbers. The calculus of time scales was initiated by S.
Hilger in his Ph.D. thesis \cite{h1} in order to unify discrete
and continuous analysis. Let $\mathbb{T}$ be a time scale.
$\mathbb{T}$ has the topology that inherits from the real numbers
with the standard topology. For $t \in \mathbb{T}$, we define the
forward jump operator $\sigma: \mathbb{T} \rightarrow \mathbb{T}$
by $\sigma(t)= \inf \left \{s \in \mathbb{T}: s >t \right\}$, and
the backward jump operator $\rho: \mathbb{T} \rightarrow
\mathbb{T}$ by $\rho(t)= \sup \left \{s \in \mathbb{T}: s < t
\right \}$.

If $\sigma(t) > t$, we say that $t$ is right-scattered, while if
$\rho(t) < t$, we say that $t$ is left-scattered. Points that are
simultaneously right-scattered and left-scattered are called
isolated. If $\sigma(t)=t$, then $t$ is called right-dense, and if
$\rho(t)=t$, then $t$ is called left-dense. Points that are
simultaneously right-dense and left-dense are called dense. Let $t
\in \mathbb{T}$, then two mappings $\mu, \nu: \mathbb{T}
\rightarrow [0, +\infty)$ are defined as follows:
$\mu(t):=\sigma(t)-t$, $\nu(t):=t-\rho(t)$.

We introduce the sets $\mathbb{T}^{k}$, $\mathbb{T}_{k}$, and
$\mathbb{T}^{k}_{k}$, which are derived from the time scale
$\mathbb{T}$, as follows. If $\mathbb{T}$ has a left-scattered
maximum $t_{1}$, then $\mathbb{T}^{k}= \mathbb{T}-\{t_{1} \}$,
otherwise $\mathbb{T}^{k}= \mathbb{T}$. If $\mathbb{T}$ has a
right-scattered minimum $t_{2}$, then $\mathbb{T}_{k}=
\mathbb{T}-\{t_{2} \}$, otherwise $\mathbb{T}_{k}= \mathbb{T}$.
Finally, we define $\mathbb{T}_{k}^{k}= \mathbb{T}^{k} \cap
\mathbb{T}_{k}$.

Throughout the text we will denote a time scales interval by
$$
[a,b]_{\mathbb{T}}=\{t\in\mathbb{T}:a\leq t \leq b\},
\quad \mbox{with}\ a,b\in\mathbb{T}.
$$

Let $f: \mathbb{T}\rightarrow \mathbb{R}$ be a real function on a
time scale $\mathbb{T}$. Then, for $t \in \mathbb{T}^{k}$ we
define $f^{\Delta}(t)$ to be the number, if one exists, such that
for all $\epsilon >0$ there is a neighborhood $U$ of $t$ such that
for all $s \in U$,
$$
\left|f(\sigma(t))-
f(s)-f^{\Delta}(t)\left(\sigma(t)-s\right)\right| \leq \epsilon
|\sigma(t)-s|.
$$
We say that $f$ is delta differentiable on $\mathbb{T}^{k}$,
provided $f^{\Delta}(t)$ exists for all $t \in \mathbb{T}^{k}$.
Similarly, for $t \in \mathbb{T}_{k}$, we define $f^{\nabla}(t)$
to be the number value, if one exists, such that for all $\epsilon
>0$, there is a neighborhood $V$ of $t$ such that for all $s \in
V$,
$$
\left|f(\rho(t))- f(s)-f^{\nabla}(t)\left(\rho(t)-s\right)\right|
\leq \epsilon |\rho(t)-s|.
$$
We say that $f$ is nabla differentiable on $\mathbb{T}_{k}$,
provided that $f^{\nabla}(t)$ exists for all $t \in
\mathbb{T}_{k}$.

Given a function $f:\mathbb{T}\rightarrow \mathbb{R}$, then we
define $f^{\sigma}: \mathbb{T}\rightarrow \mathbb{R}$ by
$f^{\sigma}(t)=f(\sigma(t))$ for all $t \in \mathbb{T}$,
\textrm{i.e.} $f^{\sigma}= f\circ \sigma$; we define $f^{\rho}:
\mathbb{T}\rightarrow \mathbb{R}$ by $f^{\rho}(t)=f(\rho(t))$ for
all $t \in \mathbb{T}$, \textrm{i.e.} $f^{\rho}= f\circ \rho$. The
following properties hold for all $t \in \mathbb{T}^{k}$:

\begin{itemize}

\item[(i)] If $f$ is delta differentiable at $t$, then $f$ is
continuous at $t$.

\item[(ii)] If $f$ is continuous at $t$ and $t$ is
right-scattered, then $f$ is delta differentiable at $t$ with
$f^{\Delta}(t)= \frac{f^{\sigma}(t)-f(t)}{\mu(t)}$.

\item[(iii)] If $f$ is right-dense, then $f$ is delta
differentiable at $t$ if and only if the limit $\lim_{s\rightarrow
t}\frac{f(t)-f(s)}{t-s}$ exists as a finite number. In this case,
$f^{\Delta}(t)= \lim_{s\rightarrow t}\frac{f(t)-f(s)}{t-s}$.

\item[(iv)] If $f$ is delta differentiable at $t$, then
$f^{\sigma}(t)= f(t)+ \mu(t)f^{\Delta}(t)$.
\end{itemize}

Similarly, given a function $f: \mathbb{T} \rightarrow
\mathbb{R}$, the following is true for all $t \in \mathbb{T}_{k}$:
\begin{itemize}

\item[(a)] If $f$ is nabla differentiable at $t$, then $f$ is
continuous at $t$.

\item[(b)] If $f$ is continuous at $t$ and $t$ is left-scattered,
then $f$ is nabla differentiable at $t$ with $f^{\nabla}(t)=
\frac{f(t)-f^{\rho}(t)}{\nu(t)}$.

\item[(c)] If $f$ is left-dense, then $f$ is nabla differentiable
at $t$ if and only if the limit $\lim_{s\rightarrow
t}\frac{f(t)-f(s)}{t-s}$ exists as a finite number. In this case,
$f^{\nabla}(t)= \lim_{s\rightarrow t}\frac{f(t)-f(s)}{t-s}.$

\item[(d)] If $f$ is nabla differentiable at $t$, then
$f^{\rho}(t)= f(t) - \nu(t)f^{\nabla}(t)$.

\end{itemize}

A function $f: \mathbb{T} \rightarrow \mathbb{R}$ is called
rd-continuous, provided it is continuous at all right-dense points
in $\mathbb{T}$ and its left-sided limits exist at all left-dense
points in $\mathbb{T}$.

A function $f: \mathbb{T} \rightarrow \mathbb{R}$ is called
ld-continuous, provided it is continuous at all left-dense points
in $\mathbb{T}$ and its right-sided limits exist finite at all
right-dense points in $\mathbb{T}$.

A function $F: \mathbb{T} \rightarrow \mathbb{R} $ is called a
delta antiderivative of $f: \mathbb{T} \rightarrow \mathbb{R}$,
provided $F^{\Delta}(t)=f(t)$ holds for all $t \in
\mathbb{T}^{k}$. Then, the delta integral of $f$ is defined by
$\int^b_a f(t)\Delta t=F(b)-F(a)$.

A function $G: \mathbb{T} \rightarrow \mathbb{R} $ is called a
nabla antiderivative of $g: \mathbb{T} \rightarrow \mathbb{R}$,
provided $G^{\nabla}(t)=g(t)$ holds for all $t \in
\mathbb{T}_{k}$. Then, the nabla integral of $g$ is defined by
$\int^b_a g(t)\nabla t=G(b)-G(a)$.

For more details on time scales we refer the reader to
\cite{a1,abra,a2,a3,a4,b1,b2}. Now, we briefly introduce the
diamond-$\alpha$ dynamic derivative and the diamond-$\alpha$
integral \cite{sfhd,Rogers}.

Let $\mathbb{T}$ be a time scale, and $t$, $s \in \mathbb{T}$.
Following \cite{Rogers}, we define $\mu_{t s} = \sigma(t)-s$,
$\eta_{t s} = \rho(t)-s$, and $f^{\diamondsuit_{\alpha}}(t)$ to be
the value, if one exists, such that for all $\epsilon >0$ there is
a neighborhood $U$ of $t$ such that for all $s \in U$
\begin{equation*}
\left| \alpha \left[f^\sigma(t) - f(s)\right] \eta_{t s} +
(1-\alpha) \left[f^\rho(t) - f(s) \right] \mu_{t s} -
f^{\diamondsuit_{\alpha}}(t) \mu_{t s} \eta_{t s} \right| <
\epsilon \left|\mu_{t s} \eta_{t s}\right| \, .
\end{equation*}
A function $f$ is said diamond-$\alpha$ differentiable on
$\mathbb{T}^k_k$ provided $f^{\diamondsuit_{\alpha}}(t)$ exists
for all $t \in \mathbb{T}^k_k$. Let $0 \leq \alpha \leq 1$. If
$f(t)$ is differentiable on $t \in \mathbb{T}^k_k$ both in the
delta and nabla senses, then $f$ is diamond-$\alpha$
differentiable at $t$ and the dynamic derivative
$f^{\diamondsuit_{\alpha}}(t)$ is given by
\begin{equation}
\label{eq:defSmp} f^{\diamondsuit_{\alpha}}(t)= \alpha
f^{\Delta}(t)+(1-\alpha)f^{\nabla}(t)
\end{equation}
(see \cite[Theorem~3.2]{Rogers}). Equality \eqref{eq:defSmp} is
the definition of $f^{\diamondsuit_{\alpha}}(t)$ found in
\cite{sfhd}. The diamond-$\alpha$ derivative reduces to the
standard $\Delta$ derivative for $\alpha =1$, or the standard
$\nabla$ derivative for $\alpha =0$. On the other hand, it
represents a ``weighted dynamic derivative'' for $\alpha \in
(0,1)$. Furthermore, the combined dynamic derivative offers a
centralized derivative formula on any uniformly discrete time
scale $\mathbb{T}$ when $\alpha=\frac{1}{2}$.

Let $f, g: \mathbb{T} \rightarrow \mathbb{R}$ be diamond-$\alpha$
differentiable at $t \in \mathbb{T}^k_k$. Then (\textrm{cf.}
\cite[Theorem~2.3]{sfhd}),

\begin{itemize}

\item[(i)] $f+g: \mathbb{T} \rightarrow \mathbb{R}$ is
diamond-$\alpha$ differentiable at $t \in \mathbb{T}^k_k$ with
$$ (f+g)^{\diamondsuit_{\alpha}}(t)=
(f)^{\diamondsuit_{\alpha}}(t)+(g)^{\diamondsuit_{\alpha}}(t);
$$

\item[(ii)] For any constant $c$, $cf: \mathbb{T} \rightarrow
\mathbb{R}$
 is diamond-$\alpha$ differentiable at $t \in \mathbb{T}^k_k$ with
$$
(cf)^{\diamondsuit_{\alpha}}(t)= c(f)^{\diamondsuit_{\alpha}}(t);
$$

\item[(iii)] $fg: \mathbb{T} \rightarrow \mathbb{R}$ is
diamond-$\alpha$ differentiable at $t \in \mathbb{T}^k_k$ with

$$
(fg)^{\diamondsuit_{\alpha}}(t)=
(f)^{\diamondsuit_{\alpha}}(t)g(t)+ \alpha
f^{\sigma}(t)(g)^{\Delta}(t) +(1-\alpha)
f^{\rho}(t)(g)^{\nabla}(t).
$$

\end{itemize}

Let $a, t \in \mathbb{T}$, and $h: \mathbb{T} \rightarrow
\mathbb{R}$. Then, the diamond-$\alpha$ integral of $h$ from $a$
to $t$ is defined by
$$
\int_{a}^{t}h(\tau) \diamondsuit_{\alpha} \tau = \alpha
\int_{a}^{t}h(\tau) \Delta \tau +(1- \alpha) \int_{a}^{t}h(\tau)
\nabla \tau, \quad 0 \leq \alpha \leq 1,
$$
provided that there exist delta and nabla integrals of $h$ on
$\mathbb{T}$. It is clear that the diamond-$\alpha$ integral of
$h$ exists when $h$ is a continuous function. We may notice that
the $\diamondsuit_{\alpha}$ combined derivative is not a dynamic
derivative for the absence of its anti-derivative
\cite[Sec.~4]{Rogers}. Moreover, in general we do not have
\begin{equation}
\label{eq:tfci:nh} \left ( \int_{a}^{t}h(\tau)
\diamondsuit_{\alpha} \tau \right)^{\diamondsuit_{\alpha}} = h(t)
\, , \quad t \in \mathbb{T} \, .
\end{equation}

\begin{example}
\label{ex:2.1}
Let $\mathbb{T} = \{0,1,2\}$, $a = 0$, and $h(\tau) = \tau^2$,
$\tau \in $ $\mathbb{T}$. It is a simple exercise to see that
$$
\left.\left ( \int_{0}^{t} h(\tau) \diamondsuit_{\alpha} \tau
\right)^{\diamondsuit_{\alpha}}\right|_{t=1} = h(1) + 2 \alpha
(1-\alpha) \, ,
$$
so that equality \eqref{eq:tfci:nh} holds only when
$\diamondsuit_{\alpha} = \nabla$ or $\diamondsuit_{\alpha} =
\Delta$.
\end{example}

Let $a$, $b$, $t \in \mathbb{T}$, $c \in \mathbb{R}$. Then
(\textrm{cf.} \cite[Theorem~3.7]{sfhd}),

\begin{itemize}

\item[(a)]$ \int_{a}^{t}\left( f(\tau)+g(\tau) \right)
\diamondsuit_{\alpha} \tau = \int_{a}^{t} f(\tau)
\diamondsuit_{\alpha} \tau + \int_{a}^{t} g(\tau)
\diamondsuit_{\alpha} \tau$;

\item[(b)] $\int_{a}^{t} c f(\tau) \diamondsuit_{\alpha} \tau = c
\int_{a}^{t} f(\tau) \diamondsuit_{\alpha} \tau$;

\item[(c)] $\int_{a}^{t} f(\tau) \diamondsuit_{\alpha} \tau =
\int_{a}^{b} f(\tau) \diamondsuit_{\alpha} \tau + \int_{b}^{t}
f(\tau) \diamondsuit_{\alpha} \tau$.

\end{itemize}

Next lemma provides some straightforward but useful results for
what follows.

\begin{lemma}\label{lem1}
Assume that $f$ and $g$ are continuous functions on
$[a,b]_{\mathbb{T}}$.
\begin{enumerate}
    \item If $f(t)\geq 0$ for all $t\in[a,b]_{\mathbb{T}}$,
    then $\int_a^b f(t)\Diamond_\alpha t\geq 0$.
    \item If $f(t)\leq g(t)$ for all $t\in[a,b]_{\mathbb{T}}$,
    then $\int_a^b f(t)\Diamond_\alpha t\leq\int_a^b g(t)\Diamond_\alpha t$.
    \item If $f(t)\geq 0$ for all $t\in[a,b]_{\mathbb{T}}$, then $f(t)=0$ if and
    only if $\int_a^b f(t)\Diamond_\alpha t=0$.
\end{enumerate}
\end{lemma}

\begin{proof}
Let $f(t)$ and $g(t)$ be continuous functions on
$[a,b]_{\mathbb{T}}$.
\begin{enumerate}
    \item Since $f(t)\geq 0$ for all $t\in[a,b]_{\mathbb{T}}$, we know (see
    \cite{b1,b2}) that $\int_a^b f(t)\Delta t\geq 0$
    and $\int_a^b f(t)\nabla t\geq 0$. Since $\alpha\in[0,1]$, the result follows.
    \item Let $h(t)=g(t)-f(t)$. Then, $\int_a^b h(t)\Diamond_\alpha t\geq
    0$ and the result follows from properties (a) and (b)
    above.
    \item If $f(t)=0$ for all $t\in[a,b]_{\mathbb{T}}$,
    the result is immediate. Suppose now that
    there exists $t_0\in[a,b]_{\mathbb{T}}$ such that $f(t_0)>0$. It is easy to
    see that at least one of the integrals $\int_a^b f(t)\Delta
    t$ or $\int_a^b f(t)\nabla t$ is strictly positive. Then,
    we have the contradiction $\int_a^b f(t)\Diamond_\alpha t>0$.
\end{enumerate}
\end{proof}

%%%%%%%%%%%%%%%%%%%%%%

\section{Main Results}
\label{sec:mr}

We now prove Jensen's diamond-$\alpha$ integral inequalities.

\begin{theorem}[Jensen's inequality]
\label{thm3} Let $\mathbb{T}$ be a time scale, $a$, $b \in
\mathbb{T}$ with $a < b$, and $c$, $d \in \mathbb{R}$. If $ g \in
C([a, b]_{\mathbb{T}}, (c, d))$ and $f \in C((c, d), \mathbb{R} )
$ is convex, then
\begin{equation}
\label{eq:JI:mr} f \left( \frac{\int_{a}^{b}g(s)
\diamondsuit_{\alpha} s}{b-a}\right ) \leq
\frac{\int_{a}^{b}f(g(s)) \diamondsuit_{\alpha} s}{b-a}.
\end{equation}
\end{theorem}

\begin{remark}
In the particular case $\alpha=1$, inequality \eqref{eq:JI:mr}
reduces to that of Theorem~\ref{thm2}. If $\mathbb{T}=\mathbb{R}$,
then Theorem~\ref{thm3} gives the classical Jensen inequality,
\textrm{i.e.} Theorem~\ref{thm1}. However, if
$\mathbb{T}=\mathbb{Z}$ and $f(x)=-\ln(x)$, then one gets the
well-known arithmetic-mean geometric-mean inequality
(\ref{eq:am:gm:i}).
\end{remark}

\begin{proof}
Since $f$ is convex we have
\begin{equation*}
\begin{split}
f \left( \frac{\int_{a}^{b}g(s) \diamondsuit_{\alpha}
s}{b-a}\right )& = f \left(\frac{\alpha}{b-a} \int_{a}^{b} g(s)
\Delta s+
\frac{1-\alpha}{b-a} \int_{a}^{b} g(s) \nabla s \right)\\
& \leq \alpha f \left(\frac{1}{b-a} \int_{a}^{b} g(s) \Delta s
\right)+ (1-\alpha)f\left( \frac{1}{b-a} \int_{a}^{b} g(s) \nabla
s \right).
\end{split}
\end{equation*}
Using now Jensen's inequality on time scales (see
Theorem~\ref{thm2}), we get
\begin{equation*}
\begin{split}
f \left( \frac{\int_{a}^{b}g(s) \diamondsuit_{\alpha}
s}{b-a}\right )& \leq \frac{\alpha}{b-a} \int_{a}^{b} f(g(s))
\Delta s +
\frac{1-\alpha}{b-a} \int_{a}^{b} f(g(s)) \nabla s \\
& = \frac{1}{b-a} \left( \alpha \int_{a}^{b} f(g(s)) \Delta s
+ (1-\alpha) \int_{a}^{b} f(g(s)) \nabla s \right) \\
& = \frac{1}{b-a} \int_{a}^{b}f(g(s)) \diamondsuit_{\alpha} s.
\end{split}
\end{equation*}
\end{proof}

Now we give an extended Jensen inequality on time scales via the
diamond-$\alpha$ integral.

\begin{theorem}[Generalized Jensen's inequality]
\label{thm4} Let $\mathbb{T}$ be a time scale, $a$, $b \in
\mathbb{T}$ with $a < b$, $c$, $d \in \mathbb{R}$, $g \in C([a,
b]_{\mathbb{T}}, (c, d))$ and $h\in C([a, b]_{\mathbb{T}},
\mathbb{R} )$ with
$$ \int_{a}^{b} |h(s)| \diamondsuit_{\alpha} s > 0 \, . $$
If $f \in C((c, d), \mathbb{R}) $ is convex, then
\begin{equation}
\label{eq:GJI} f \left( \frac{\int_{a}^{b} |h(s)|g(s)
\diamondsuit_{\alpha} s}{\int_{a}^{b} |h(s)|
\diamondsuit_{\alpha}s}\right ) \leq \frac{\int_{a}^{b}
|h(s)|f(g(s)) \diamondsuit_{\alpha} s}{\int_{a}^{b} |h(s)|
\diamondsuit_{\alpha}s}.
\end{equation}
\end{theorem}

\begin{remark}
Theorem~\ref{thm4} is the same as
\cite[Theorem~3.17]{rev1:r}. However,
we prove Theorem~\ref{thm4} using a different approach
than that proposed in \cite{rev1:r}:
in \cite{rev1:r} it is stated
that such result follows from the analog nabla-inequality.
As we have seen, diamond-alpha integrals have different
properties than those of delta or nabla integrals
(\textrm{cf.} Example~\ref{ex:2.1}). On the other hand,
there is an inconsistency in \cite{rev1:r}: a very simple example showing this fact is given below in Remark~\ref{rem:ex:er}.
\end{remark}

\begin{remark}
In the particular case $h=1$, Theorem~\ref{thm4} reduces to
Theorem~\ref{thm3}.
\end{remark}

\begin{remark}
If $f$ is strictly convex, the inequality sign ``$\leq$'' in
(\ref{eq:GJI}) can be replaced by ``$<$''. Similar result to
Theorem~\ref{thm4} holds if one changes the condition ``$f$ is
convex'' to ``$f$ is concave'', by replacing the inequality sign
``$\leq$'' in (\ref{eq:GJI}) by ``$\geq$''.
\end{remark}

\begin{proof}
Since $f$ is convex, it follows, for example from \cite[Exercise
3.42C]{fol}, that for $t \in (c, d)$ there exists $a_{t} \in
\mathbb{R} $ such that
\begin{equation}\label{eq1}
a_{t}(x-t) \leq f(x)-f(t) \mbox{ for all } x \in (c, d).
\end{equation}
Setting $$ t= \frac{\int_{a}^{b} |h(s)|g(s) \diamondsuit_{\alpha}
s}{\int_{a}^{b} |h(s)| \diamondsuit_{\alpha}s} \, ,
$$
then using \eqref{eq1} and item 2 of Lemma~\ref{lem1}, we get
\begin{equation*}
\begin{split}
\int_{a}^{b} & |h(s)|f(g(s)) \diamondsuit_{\alpha} s - \left (
\int_{a}^{b} |h(s)| \diamondsuit_{\alpha} s \right) f \left(
\frac{\int_{a}^{b} |h(s)|g(s) \diamondsuit_{\alpha}
s}{\int_{a}^{b}
|h(s)| \diamondsuit_{\alpha}s}\right ) \\
 = & \int_{a}^{b} |h(s)|f(g(s)) \diamondsuit_{\alpha} s - \left (
\int_{a}^{b} |h(s)| \diamondsuit_{\alpha} s \right) f(t)
=  \int_{a}^{b} |h(s)| \left (f(g(s))-f(t) \right)
\diamondsuit_{\alpha} s \\
\geq & a_{t} \left( \int_{a}^{b} |h(s)| (g(s)-t ) \right)
\diamondsuit_{\alpha} s
=  a_{t} \left( \int_{a}^{b} |h(s)|g(s) \diamondsuit_{\alpha}
s - t \int_{a}^{b} |h(s)| \diamondsuit_{\alpha} s \right)\\
\\ = & a_{t} \left( \int_{a}^{b} |h(s)|g(s) \diamondsuit_{\alpha}
s- \int_{a}^{b} |h(s)|g(s) \diamondsuit_{\alpha} s \right) = 0 \,
.
\end{split}
\end{equation*}
This leads to the desired inequality.
\end{proof}

\begin{remark}
The proof of Theorem~\ref{thm4} follows closely the proof of the
classical Jensen inequality (see \textrm{e.g.} \cite[Problem
3.42]{fol}) and the proof of Jensen's inequality on time scales
\cite{abp}.
\end{remark}

We have the following corollaries.

\begin{corollary}$(\mathbb{T}=\mathbb{R})$
Let $g, h: [a, b]\rightarrow \mathbb{R}$ be continuous functions
with $g([a, b]) \subseteq (c, d)$ and $\int_{a}^{b}|h(x)| dx >0$.
If $f \in C((c, d), \mathbb{R})$ is convex, then
$$
f \left( \frac{\int_{a}^{b} |h(x)|g(x) dx}{\int_{a}^{b} |h(x)|
dx}\right ) \leq \frac{\int_{a}^{b} |h(x)|f(g(x)) dx}{\int_{a}^{b}
|h(x)| dx}.
$$
\end{corollary}

\begin{corollary}$(\mathbb{T}=\mathbb{Z})$
\label{cor:E:SMC} Given a convex function $f$, we have for any
$x_{1}, \ldots ,x_{n} \in \mathbb{R}$ and $c_{1}, \ldots ,c_{n}
\in \mathbb{R}$ with $\sum_{k=1}^{n}|c_{k}|
>0$:
\begin{equation}
\label{eq:des:pr} f \left
(\frac{\sum_{k=1}^{n}|c_{k}|x_{k}}{\sum_{k=1}^{n}|c_{k}|}\right)
\leq \frac{\sum_{k=1}^{n}|c_{k}|f(x_{k})}{\sum_{k=1}^{n}|c_{k}|}.
\end{equation}
\end{corollary}

\begin{remark}
\label{rem:ex:er}
Corollary~\ref{cor:E:SMC} coincides with \cite[Corollary~2.4]{fcw}
and \cite[Corollary~3.12]{rev1:r} if one substitutes all the
$|c_{k}|$'s in Corollary~\ref{cor:E:SMC} by $c_k$ and we restrict
ourselves to integer values of $x_i$ and $c_i$, $i = 1,\ldots,n$.
Let $\mathbb{T}=\mathbb{Z}$, $a = 1$ and $b = 3$, so that
$[a,b]_{\mathbb{T}}$ denotes the set $\{1,2,3\}$ and $n=3$.
For the data $f(x)=x^2$, $c_1=1$, $c_2=5$, $c_3=-3$, $x_1=1$, $x_2=1$, and
$x_3=2$ one has $A=\sum_{k=1}^3 c_k = 3 > 0$ and $B=\sum_{k=1}^3
c_k x_k = 0$. Thus, $D=f(B/A)=f(0)=0$. On the other hand,
$f(x_1)=1$, $f(x_2)=1$, and $f(x_3)=4$. Therefore, $C=\sum_{k=1}^3
c_k f(x_k)= -6$. We have $E=C/A=-2$ and $D>E$, \textrm{i.e.} $f
\left (\frac{\sum_{k=1}^{n} c_{k} x_{k}}{\sum_{k=1}^{n}
c_{k}}\right) > \frac{\sum_{k=1}^{n} c_{k}
f(x_{k})}{\sum_{k=1}^{n} c_{k}}$. Inequality \eqref{eq:des:pr}
gives the truism $\frac{16}{9} \le 2$.
\end{remark}

%%%%%%%%%%%%%%%%%%%%%%

\section*{Particular Cases}

\begin{itemize}
\item[(i)] Let $g(t) > 0$ on $[a, b]_{\mathbb{T}}$ and $f(t)=
t^{\beta}$ on $(0, +\infty)$. One can see that $f$ is convex on
$(0, +\infty)$ for $\beta < 0$ or $\beta >1$, and $f$ is concave
on $(0, +\infty)$ for $\beta \in (0, 1)$. Then,
$$
 \left( \frac{\int_{a}^{b} |h(s)|g(s) \diamondsuit_{\alpha}
s}{\int_{a}^{b} |h(s)| \diamondsuit_{\alpha}s}\right )^{\beta}
\leq \frac{\int_{a}^{b} |h(s)|g^{\beta}(s) \diamondsuit_{\alpha}
s}{\int_{a}^{b} |h(s)| \diamondsuit_{\alpha}s}, \mbox{ if } \beta
< 0 \mbox{ or } \beta >1;
$$

$$
 \left( \frac{\int_{a}^{b} |h(s)|g(s) \diamondsuit_{\alpha}
s}{\int_{a}^{b} |h(s)| \diamondsuit_{\alpha}s}\right )^{\beta}
\geq \frac{\int_{a}^{b} |h(s)|g^{\beta}(s) \diamondsuit_{\alpha}
s}{\int_{a}^{b} |h(s)| \diamondsuit_{\alpha}s}, \mbox{ if } \beta
\in (0, 1).
$$

\item[(ii)] Let $g(t) > 0$ on $[a, b]_{\mathbb{T}}$ and $f(t)=
\ln(t)$ on $(0, +\infty)$. One can also see that $f$ is concave on
$(0, +\infty)$. It follows that

$$
 \ln \left( \frac{\int_{a}^{b} |h(s)|g(s) \diamondsuit_{\alpha}
s}{\int_{a}^{b} |h(s)| \diamondsuit_{\alpha}s}\right ) \geq
\frac{\int_{a}^{b} |h(s)|\ln (g(s)) \diamondsuit_{\alpha}
s}{\int_{a}^{b} |h(s)| \diamondsuit_{\alpha}s}.
$$

\item[(iii)]Let $h=1$, then
$$
 \ln \left(\frac{\int_{a}^{b} g(s) \diamondsuit_{\alpha}
s}{b - a}\right) \geq \frac{\int_{a}^{b} \ln (g(s)) \diamondsuit_{\alpha} s}{b - a}
.
$$

\item[(iv)] Let $\mathbb{T}=\mathbb{R}$, $g: [0, 1]\rightarrow (0,
\infty)$ and $h(t)=1$. Applying Theorem~\ref{thm4} with the convex
and continuous function $f=-\ln$ on $(0, \infty)$, $a=0$ and
$b=1$, we get:
$$
 \ln \int_{0}^{1} g(s)ds \geq \int_{0}^{1}\ln( g(s))ds.
$$
Then,
$$
\int_{0}^{1} g(s)ds \geq \exp \left(\int_{0}^{1}\ln( g(s))ds
\right).
$$
\item[(v)] Let $\mathbb{T}=\mathbb{Z}$ and $n\in\mathbb{N}$. Fix
$a=1$, $b=n+1$ and consider a function
$g:\{1,\ldots,n+1\}\rightarrow(0,\infty)$. Obviously, $f=-\ln$ is
convex and continuous on $(0,\infty)$, so we may apply Jensen's
inequality to obtain
\end{itemize}
\begin{equation*}
\begin{split}
\ln\Biggl[ & \frac{1}{n}\left(\alpha\sum_{t=1}^n
g(t)+(1-\alpha)\sum_{t=2}^{n+1}g(t)\right)\Biggr] =
\ln\left[\frac{1}{n}\int_1^{n+1}g(t)\Diamond_\alpha
t\right]\\
&\geq\frac{1}{n}\int_1^{n+1}\ln(g(t))\Diamond_\alpha t\\
&=\frac{1}{n}\left[\alpha\sum_{t=1}^n
\ln(g(t))+(1-\alpha)\sum_{t=2}^{n+1}\ln(g(t))\right]\\
&=\ln\left\{\prod_{t=1}^n
g(t)\right\}^{\frac{\alpha}{n}}+\ln\left\{\prod_{t=2}^{n+1}
g(t)\right\}^{\frac{1-\alpha}{n}} \, ,
\end{split}
\end{equation*}
and hence
$$\frac{1}{n}\left(\alpha\sum_{t=1}^n
g(t)+(1-\alpha)\sum_{t=2}^{n+1}g(t)\right)\geq\left\{\prod_{t=1}^n
g(t)\right\}^{\frac{\alpha}{n}}\left\{\prod_{t=2}^{n+1}
g(t)\right\}^{\frac{1-\alpha}{n}}.$$ When $\alpha=1$, we obtain
the well-known arithmetic-mean geometric-mean inequality:
\begin{equation}
\label{eq:am:gm:i} \frac{1}{n}\sum_{t=1}^n
g(t)\geq\left\{\prod_{t=1}^n g(t)\right\}^{\frac{1}{n}}.
\end{equation}
When $\alpha=0$, we also have $$\frac{1}{n}\sum_{t=2}^{n+1}
g(t)\geq\left\{\prod_{t=2}^{n+1} g(t)\right\}^{\frac{1}{n}}.$$

\begin{itemize}

\item[(vi)] Let $\mathbb{T}= 2^{\mathbb{N}_{0}}$ and $N \in
\mathbb{N}$. We can apply Theorem~\ref{thm4} with $a=1, b=2^{N}$
and $g: \{ 2^{k}: 0 \leq k \leq N \} \rightarrow (0, \infty)$.
Then, we get:
\end{itemize}

\begin{equation*}
\begin{split}
\ln \left \{
\frac{\int_{1}^{2^{N}}g(t)\diamondsuit_{\alpha}t}{2^{N}-1} \right
\}&= \ln \left \{ \alpha \frac{\int_{1}^{2^{N}}g(t)\Delta
t}{2^{N}-1}+(1-\alpha) \frac{\int_{1}^{2^{N}}g(t)\nabla
t}{2^{N}-1}
\right \}\\
&= \ln \left \{ \frac{\alpha
\sum_{k=0}^{N-1}2^{k}g(2^{k})}{2^{N}-1}
+ \frac{(1-\alpha) \sum_{k=1}^{N}2^{k}g(2^{k})}{2^{N}-1} \right \}\\
& \geq \frac{\int_{1}^{2^{N}} \ln
(g(t))\diamondsuit_{\alpha}t}{2^{N}-1}
\end{split}
\end{equation*}

\begin{equation*}
\begin{split}
&= \alpha \frac{\int_{1}^{2^{N}} \ln (g(t))\Delta t}{2^{N}-1}+
(1 - \alpha) \frac{\int_{1}^{2^{N}} \ln (g(t))\nabla t}{2^{N}-1}\\
&= \frac{\alpha \sum_{k=0}^{N-1}2^{k}\ln(g(2^{k})) }{2^{N}-1} +
\frac{(1-\alpha) \sum_{k=1}^{N}2^{k}\ln(g(2^{k})) }{2^{N}-1}\\
&= \frac{ \sum_{k=0}^{N-1}\ln(g(2^{k}))^{\alpha 2^{k}} }{2^{N}-1}
+
\frac{ \sum_{k=1}^{N}\ln(g(2^{k}))^{(1-\alpha)2^{k}} }{2^{N}-1}\\
&=\frac{ \ln \prod_{k=0}^{N-1}(g(2^{k}))^{\alpha 2^{k}} }{2^{N}-1}
+ \frac{ \ln(\prod_{k=1}^{N}g(2^{k}))^{(1-\alpha)2^{k}}
}{2^{N}-1}\\
&= \ln \left \{\prod_{k=0}^{N-1}(g(2^{k}))^{\alpha 2^{k}} \right
\}^{\frac{1}{2^{N}-1}} + \ln \left
\{\prod_{k=1}^{N}(g(2^{k}))^{(1-\alpha) 2^{k}} \right
\}^{\frac{1}{2^{N}-1}}\\
&= \ln \left ( \left \{\prod_{k=0}^{N-1}(g(2^{k}))^{\alpha 2^{k}}
\right \}^{\frac{1}{2^{N}-1}} \left
\{\prod_{k=1}^{N}(g(2^{k}))^{(1-\alpha) 2^{k}} \right
\}^{\frac{1}{2^{N}-1}} \right) \, .
\end{split}
\end{equation*}
We conclude that
\begin{multline*}
  \ln \left \{ \frac{\alpha
\sum_{k=0}^{N-1}2^{k}g(2^{k})+(1-\alpha)
\sum_{k=1}^{N}2^{k}g(2^{k})}{2^{N}-1} \right \}\\
\geq \ln \left ( \left \{\prod_{k=0}^{N-1}(g(2^{k}))^{\alpha
2^{k}} \right \}^{\frac{1}{2^{N}-1}} \left
\{\prod_{k=1}^{N}(g(2^{k}))^{(1-\alpha) 2^{k}} \right
\}^{\frac{1}{2^{N}-1}} \right).
\end{multline*}
On the other hand,
$$
\alpha \sum_{k=0}^{N-1}2^{k}g(2^{k})+(1-\alpha)
\sum_{k=1}^{N}2^{k}g(2^{k})= \sum_{k=1}^{N-1}2^{k}g(2^{k})+\alpha
g(1)+(1-\alpha) 2^{N}g(2^{N}).
$$
It follows that
\begin{multline*}
\frac{\sum_{k=1}^{N-1}2^{k}g(2^{k})+\alpha g(1)+(1-\alpha)
2^{N}g(2^{N})}{2^{N}-1} \\
\geq \left\{\prod_{k=0}^{N-1}(g(2^{k}))^{\alpha 2^{k}} \right
\}^{\frac{1}{2^{N}-1}} \left
\{\prod_{k=1}^{N}(g(2^{k}))^{(1-\alpha) 2^{k}} \right
\}^{\frac{1}{2^{N}-1}}.
\end{multline*}
In the particular case when $\alpha =1$ we have
$$
\frac{\sum_{k=0}^{N-1}2^{k}g(2^{k})}{2^{N}-1} \geq \left
\{\prod_{k=0}^{N-1}(g(2^{k}))^{ 2^{k}} \right
\}^{\frac{1}{2^{N}-1}},
$$
and when $\alpha =0$ we get the inequality
$$
\frac{\sum_{k=1}^{N}2^{k}g(2^{k})}{2^{N}-1} \geq \left
\{\prod_{k=1}^{N}(g(2^{k}))^{ 2^{k}} \right
\}^{\frac{1}{2^{N}-1}}\, .
$$

%%%%%%%%%%%%%%%%%%%%%%

\section{Related Diamond-$\alpha$ Integral Inequalities}
\label{sec:app}

The usual proof of H\"{o}lder's inequality use the basic Young
inequality $x^{\frac{1}{p}}y^{\frac{1}{q}} \leq \frac{x}{p}+
\frac{y}{q}$ for nonnegative $x$ and $y$. Here we present a proof
based on the application of Jensen's inequality
(Theorem~\ref{thm4}).

\begin{theorem}[H\"{o}lder's inequality]
\label{app:th:hi} Let $\mathbb{T}$ be a time scale, $a$, $b \in
\mathbb{T}$ with $a < b$, and $f$, $g$, $h \in C([a,
b]_{\mathbb{T}}, [0, \infty))$ with
$\int_{a}^{b}h(x)g^{q}(x)\diamondsuit_{\alpha} x
>0$, where $q$ is the H\"{o}lder conjugate number of $p$,
\textrm{i.e.} $\frac{1}{p}+\frac{1}{q}=1$ with $1<p$. Then, we
have:
\begin{equation}
\label{app:eq:hi} \int_{a}^{b}h(x)f(x)g(x)\diamondsuit_{\alpha} x
\leq \left(\int_{a}^{b}h(x)f^{p}(x)\diamondsuit_{\alpha}
x\right)^{\frac{1}{p}}
\left(\int_{a}^{b}h(x)g^{q}(x)\diamondsuit_{\alpha}
x\right)^{\frac{1}{q}} \, .
\end{equation}
\end{theorem}

\begin{proof}
Choosing $f(x)=x^{p}$ in Theorem~\ref{thm4}, which for $p>1$ is
obviously a convex function on $[0, \infty)$, we have
\begin{equation}
\label{eq:R} \left( \frac{\int_{a}^{b} |h(s)|g(s)
\diamondsuit_{\alpha} s}{\int_{a}^{b} |h(s)|
\diamondsuit_{\alpha}s}\right )^p \leq \frac{\int_{a}^{b}
|h(s)|(g(s))^p \diamondsuit_{\alpha} s}{\int_{a}^{b} |h(s)|
\diamondsuit_{\alpha}s}.
\end{equation}
Inequality \eqref{app:eq:hi} is trivially true in the case when
$g$ is identically zero. We consider two cases: (i) $g(x) > 0$ for
all $x \in [a, b]_{\mathbb{T}}$; (ii) there exists at least one $x
\in [a, b]_{\mathbb{T}}$ such that $g(x) = 0$. We begin with
situation (i). Replacing $g$ by $fg^{\frac{-q}{p}}$ and $|h(x)|$
by $hg^{q}$ in inequality (\ref{eq:R}), we get:
$$
\left(
\frac{\int_{a}^{b}h(x)g^{q}(x)f(x)g^{\frac{-q}{p}}(x)\diamondsuit_{\alpha}
x}{\int_{a}^{b}h(x)g^{q}(x)\diamondsuit_{\alpha}x} \right)^{p}
\leq
\frac{\int_{a}^{b}h(x)g^{q}(x)(f(x)g^{\frac{-q}{p}}(x))^{p}\diamondsuit_{\alpha}
x}{\int_{a}^{b}h(x)g^{q}(x)\diamondsuit_{\alpha}x}.
$$
Using the fact that $\frac{1}{p}+\frac{1}{q}=1$, we obtain that
\begin{equation}
\label{eq:parti} \int_{a}^{b}h(x)f(x)g(x)\diamondsuit_{\alpha} x
\leq \left(\int_{a}^{b}h(x)f^{p}(x)\diamondsuit_{\alpha}
x\right)^{\frac{1}{p}}
\left(\int_{a}^{b}h(x)g^{q}(x)\diamondsuit_{\alpha}
x\right)^{\frac{1}{q}} \, .
\end{equation}
We now consider situation (ii). Let $G= \left\{x \in [a,
b]_{\mathbb{T}} \, | \, g(x) =0 \right\}$. Then,
\begin{gather*}
\int_{a}^{b}h(x) f(x) g(x) \diamondsuit_{\alpha} x = \int_{[a,
b]_{\mathbb{T}}-G} h(x) f(x) g(x) \diamondsuit_{\alpha} x +
\int_{G} h(x) f(x) g(x) \diamondsuit_{\alpha} x\\
= \int_{[a, b]_{\mathbb{T}}-G} h(x) f(x) g(x)
\diamondsuit_{\alpha} x
\end{gather*}
because $\int_{G} h(x) f(x) g(x) \diamondsuit_{\alpha} x =0$. For the set $[a, b]_{\mathbb{T}}-G$ we are in case (i), \textrm{i.e.}
$g(x) > 0$, and it follows from \eqref{eq:parti} that
\begin{equation*}
\begin{split}
\int_{a}^{b}h(x) f(x) g(x) \diamondsuit_{\alpha} x &= \int_{[a,
b]_{\mathbb{T}}-G}h(x) f(x) g(x) \diamondsuit_{\alpha} x \\
&\leq \left(\int_{[a, b]_{\mathbb{T}}-G} h(x) f^{p}(x)
\diamondsuit_{\alpha} x \right )^{\frac{1}{p}} \quad \left
(\int_{[a, b]_{\mathbb{T}}-G} h(x)
g^{q}(x) \diamondsuit_{\alpha} x \right )^{\frac{1}{q}}\\
&\leq \left (\int_a^b h(x) f^{p}(x) \diamondsuit_{\alpha} x \right
)^{\frac{1}{p}} \quad \left (\int_a^b h(x) g^{q}(x)
\diamondsuit_{\alpha} x \right )^{\frac{1}{q}} \, .
\end{split}
\end{equation*}
\end{proof}

\begin{remark}
In the particular case $h=1$, Theorem~\ref{app:th:hi} gives the
diamond-$\alpha$ version of classical H\"{o}lder's inequality:
$$
\int_{a}^{b}|f(x)g(x)|\diamondsuit_{\alpha} x \leq
\left(\int_{a}^{b}|f|^{p}(x)\diamondsuit_{\alpha}
x\right)^{\frac{1}{p}}
\left(\int_{a}^{b}|g|^{q}(x)\diamondsuit_{\alpha}
x\right)^{\frac{1}{q}},
$$
where $p>1$ and $q=\frac{p}{p-1}$.
\end{remark}

\begin{remark}
In the special case $p=q=2$, (\ref{app:eq:hi}) reduces to the
following diamond-$\alpha$ Cauchy-Schwarz integral inequality on
time scales:
$$
\int_{a}^{b}|f(x)g(x)|\diamondsuit_{\alpha} x \leq \sqrt{
\left(\int_{a}^{b}f^{2}(x)\diamondsuit_{\alpha} x\right)
 \left(\int_{a}^{b}g^{2}(x)\diamondsuit_{\alpha}
x\ \right)} \, .
$$
\end{remark}

We are now in position to prove a Minkowski inequality using our
H\"{o}lder's inequality (\ref{app:eq:hi}).

\begin{theorem}[Minkowski's inequality]
Let $\mathbb{T}$ be a time scale, $a$, $b \in \mathbb{T}$ with $a
< b$, and $p>1$. For continuous functions $f, g: [a,
b]_{\mathbb{T}} \rightarrow \mathbb{R}$ we have
$$
\left (\int_{a}^{b}|(f+g)(x)|^{p}\diamondsuit_{\alpha} x
\right)^{\frac{1}{p}} \leq
\left(\int_{a}^{b}|f(x)|^{p}\diamondsuit_{\alpha}
x\right)^{\frac{1}{p}} +
\left(\int_{a}^{b}|g(x)|^{p}\diamondsuit_{\alpha}
x\right)^{\frac{1}{p}}.
$$
\end{theorem}

\begin{proof}
We have, by the triangle inequality, that
\begin{multline}
\label{mink1} \int_a^b |f(x) + g(x)|^p\diamondsuit_{\alpha} x
=\int_a^b |f(x)+g(x)|^{p-1}|f(x)+g(x)|\diamondsuit_\alpha x\\
\leq \int_a^b|f(x)||f(x)+g(x)|^{p-1}\diamondsuit_\alpha
x+\int_a^b|g(x)||f(x)+g(x)|^{p-1}\diamondsuit_\alpha x.
\end{multline}
Applying now H\"{o}lder's inequality with $q=p/(p-1)$ to
\eqref{mink1}, we obtain:
\begin{multline*}
\int_a^b |f(x)+g(x)|^p\diamondsuit_{\alpha} x
\leq\left[\int_a^b|f(x)|^p\diamondsuit_\alpha
x\right]^{\frac{1}{p}}\left[\int_a^b
|f(x)+g(x)|^{(p-1)q}\diamondsuit_\alpha
x\right]^{\frac{1}{q}}\\
 +\left[\int_a^b|g(x)|^p\diamondsuit_\alpha
x\right]^{\frac{1}{p}}\left[\int_a^b
|f(x)+g(x)|^{(p-1)q}\diamondsuit_\alpha
x\right]^{\frac{1}{q}}\\
=\left\{\left[\int_a^b|f(x)|^p\diamondsuit_\alpha
x\right]^{\frac{1}{p}}+\left[\int_a^b|g(x)|^p\diamondsuit_\alpha
x\right]^{\frac{1}{p}}\right\}\left[\int_a^b|f(x)+g(x)|^p\diamondsuit_\alpha
x\right]^{\frac{1}{q}}.
\end{multline*}
Dividing both sides of the last inequality by
$$\left[\int_a^b|f(x)+g(x)|^p\diamondsuit_\alpha
x\right]^{\frac{1}{q}},$$ we get the desired conclusion.
\end{proof}

As another application of Theorem~\ref{thm4}, we have:

\begin{theorem}
\label{thm5} Let $\mathbb{T}$ be a time scale, $a$, $b \in
\mathbb{T}$ with $a < b$, and $f$, $g$, $h \in C([a,
b]_{\mathbb{T}}, [0, \infty))$.
\begin{itemize}
\item[(i)] If $p>1$, then
$$
\left\{\left(\int_{a}^{b} hf \diamondsuit_{\alpha}x\right)^{p}
+\left(\int_{a}^{b} hg \diamondsuit_{\alpha}x\right)^{p}
\right\}^{\frac{1}{p}} \leq \int_{a}^{b}
h(f^{p}+g^{p})^{\frac{1}{p}}\diamondsuit_{\alpha}x \, .
$$
\item[(ii)] If $\ 0 < p < 1$, then
$$ \left\{\left(\int_{a}^{b} hf
\diamondsuit_{\alpha}x\right)^{p} +\left(\int_{a}^{b} hg
\diamondsuit_{\alpha}x\right)^{p} \right\}^{\frac{1}{p}} \geq
\int_{a}^{b} h(f^{p}+g^{p})^{\frac{1}{p}}\diamondsuit_{\alpha}x \,
.
$$
\end{itemize}
\end{theorem}

\begin{proof}
We prove only (i). The proof of (ii) is similar. Inequality (i) is trivially true when $f$ is zero: both the left and
right hand sides reduce to $\int_a^b h g \diamondsuit_{\alpha}x$.
Otherwise, applying Theorem~\ref{thm4} with
$f(x)=(1+x^{p})^{\frac{1}{p}}$, which is clearly convex on $(0,
\infty)$, we obtain
$$
\left(1+\frac{(\int_{a}^{b}hf\diamondsuit_{\alpha}x)^{p}}
{(\int_{a}^{b}h\diamondsuit_{\alpha}x)^{p}}\right)^{\frac{1}{p}}
\leq
\frac{\int_{a}^{b}h(1+f^{p})^{\frac{1}{p}}\diamondsuit_{\alpha}x}
{\int_{a}^{b}h\diamondsuit_{\alpha}x}.
$$
In other words,
$$
\left[ \left(\int_{a}^{b}h\diamondsuit_{\alpha}x\right)^{p}+
\left(\int_{a}^{b}hf\diamondsuit_{\alpha}x\right)^{p}
\right]^{\frac{1}{p}} \leq
\int_{a}^{b}h(1+f^{p})^{\frac{1}{p}}\diamondsuit_{\alpha}x.
$$
Changing $h$ and $f$ by $\frac{hf}{\int_{a}^{b}hf
\diamondsuit_{\alpha}x}$ and $\frac{g}{f}$ in the last inequality,
respectively, we obtain directly the inequality (i) of
Theorem~\ref{thm5}.
\end{proof}

%%%%%%%%%%%%%%%%%%%%%

\section*{Acknowledgements}

The authors were supported by the \emph{Portuguese Foundation for
Science and Technology} (FCT), through the \emph{Centre for
Research on Optimization and Control} (CEOC) of the University of
Aveiro, cofinanced by the European Community Fund FEDER/POCI 2010
(all the three authors); the postdoc fellowship
SFRH/BPD/20934/2004 (Sidi Ammi); the PhD fellowship
SFRH/BD/39816/2007 (Ferreira); and the research project
PTDC/MAT/72840/2006 (Torres). The authors are grateful to three referees, and the editor assigned to handle the review process, for several helpful comments and a careful reading of the manuscript.

%%%%%%%%%%%%%%%%%%%%%

%%%%%%%%%%%%%%%%%%%%%

\end{document}